\magnification = 1150
\baselineskip = 14 pt

\def \sect#1{\bigskip \noindent {\bf #1} \medskip}
\def \subsect#1{\medskip \noindent{\it #1} \medskip}
\def \th#1#2{\medskip \noindent {\bf Theorem #1.}   \it #2 \rm}
\def \prop#1#2{\medskip \noindent {\bf Proposition #1.}   \it #2 \rm}
\def \cor#1#2{\medskip \noindent {\bf Corollary #1.}   \it #2 \rm}
\def \pf {\noindent  {\it Proof}.\quad }
\def \lem#1#2{\medskip \noindent {\bf Lemma #1.}   \it #2 \rm}

\def\sqr#1#2{{\vcenter{\vbox{\hrule height.#2pt\hbox{\vrule width.#2pt height#1pt \kern#1pt\vrule width.#2pt}\hrule height.#2pt}}}}

\def \square{\hfill\mathchoice\sqr56\sqr56\sqr{4.1}5\sqr{3.5}5}

\nopagenumbers

\headline={\ifnum\pageno=1 \hfill \else \hfill {\rm \folio} \fi}

\centerline{\bf Correspondence between Lifetime Minimum Wealth}
\centerline{\bf and Utility of Consumption} \bigskip

\noindent Erhan Bayraktar \hfill \break
\indent Department of Mathematics \hfill \break
\indent University of Michigan \hfill \break
\indent Ann Arbor, Michigan, 48109 \hfill \break
\indent erhan@umich.edu

\medskip

\noindent Virginia R. Young \hfill \break
\indent Department of Mathematics \hfill \break
\indent University of Michigan \hfill \break
\indent Ann Arbor, Michigan, 48109 \hfill \break
\indent vryoung@umich.edu

\bigskip

\centerline{Version: 28 November 2006} \bigskip

\noindent{\bf Abstract:}  We establish when the two problems of minimizing a function of lifetime minimum wealth and of maximizing utility of lifetime consumption result in the same optimal investment strategy on a given open interval $O$ in wealth space.  To answer this question, we equate the two investment strategies and show that if the individual consumes at the same rate in both problems -- the consumption rate is a control in the problem of maximizing utility -- then the investment strategies are equal only when the consumption function is linear in wealth on $O$, a rather surprising result.  It, then, follows that  the corresponding investment strategy is also linear in wealth and the implied utility function exhibits hyperbolic absolute risk aversion.

\medskip


\medskip


\medskip


\sect{1. Introduction}

When an individual seeks to find an optimal investment policy, the resulting optimal policy depends on her optimization criterion.  The most common optimization criterion encountered in the finance literature is to  maximize one's expected discounted utility of consumption and bequest.  Merton (1992) studies this problem, and many others have continued his work; see, for example, Karatzas and Shreve (1998, Chapter 3) and the discussion at the end of that chapter for further references.  More recently, researchers have begun to find the optimal investment policy to minimize the probability that an individual runs out of money before dying, also called the problem of minimizing the probability of lifetime ruin; see, for example, Milevsky, Ho, and Robinson (1997), Milevsky and Robinson (2000), Young (2004), and Milevsky, Moore, and Young (2006).  Note that whether someone ruins is a function of minimum wealth.

We consider the following two problems in this paper:

\noindent \hangindent 20 pt {\bf Problem 1}  An investor adopts a rate of consumption $c$ that depends on wealth and then considers the problem of minimizing a nonnegative, nonincreasing function $f$ of lifetime minimum wealth.  The consumption rate is an exogenous variable in this problem, and the investor only determines an optimal investment strategy $\pi^*$.

\noindent \hangindent 20 pt {\bf Problem 2}  An investor considers the classical Merton problem that aims at maximizing the expected discounted utility derived by the investor's rate of lifetime consumption.  Given a utility function $u$, the investor then determines an optimal investment strategy $\pi^u$ and consumption rate $c^u$.

\medskip

We establish when these two investors behave similarly (on some open interval $O$ of wealth space).  We are motivated to do this because Young (2004) and Bayraktar and Young (2006) find such a correspondence in some special cases.  In both papers, the individual seeks to minimize the probability that wealth reaches some point $b > 0$ before she dies.  Young (2004) places no restriction on the optimal investment strategy, while Bayraktar and Young (2006) extend Young's work to two cases:  (1)  The individual may not borrow any money; and (2) the individual may borrow money but only at a rate higher than the rate earned by an investment in the riskless asset.  Fleming and Zariphopoulou (1991) consider the latter setting in the problem of maximizing expected utility of consumption under power utility.  In each case, when the consumption rate is proportional to wealth, the individual who minimizes her probability of lifetime ruin behaves like an individual who maximizes her expected discounted utility of consumption when utility is a power function.

Motivated by these results, in this paper, we {\it completely} characterize the utility functions (and corresponding consumption functions) for which an investor behaves the same under Problem 1 and Problem 2.  To this end, we begin in Section 2 by solving Problem 1 for a nonincreasing, nonnegative function of lifetime minimum wealth in a Black-Scholes market.  In Section 3, we examine Problem 2 and show that if $\pi^* = \pi^u$ and if $c = c^u$ on $O$, then the derivative $c'$ is constant, a rather surprising result.  By using this result, in Section 4.1, we show that $\pi^*$ is also linear in wealth by reconsidering Problem 1 for the special case in which $c$ is linear in wealth on $O$.  In Section 4.2, by using the linearity of $\pi^u$ and $c$ in wealth, we show that the utility function $u$ exhibits hyperbolic absolute risk aversion (HARA) on $O$.  Section 5 concludes the paper.

Our main contribution lies in the result that investors under Problems 1 and 2 behave similarly if and only if the utility function is HARA (for Problem 2), summarized in Theorem 4.3.  Arguably another major contribution lies within the solution of Problem 1 in Section 2.  Specifically, we proved that the optimal investment strategy $\pi^*$ is independent of $f$, another surprising result; see Theorem 2.5.  Its derivation avoids the usual verification argument for $L^{\infty}$ control problems that exhibit a Neumann condition; see Section 2.2.

\sect{2. Lifetime Minimum Wealth in a Black-Scholes Market}

In this section, we consider the problem of minimizing the expectation of a nonincreasing, nonnegative function of lifetime minimum wealth.  In Section 2.1, we present the financial market and define lifetime minimum wealth.  In this section, we also indicate how the value function of this problem is related to the minimum probability of ruin.  In Section 2.2, we present a verification lemma that specifies a Hamilton-Jacobi-Bellman (HJB) equation that the probability of ruin satisfies, namely Theorem 2.1.  Corollary 2.6 shows that we can represent the value function of the more general problem in terms of the probability of ruin.  Later, in Section 4.1, we show how to calculate the probability of ruin for a specific rate of consumption.

\subsect{2.1. Financial Market and Definition of the Value Function $V^f$}

In this section, we present the financial ingredients that affect the individual's wealth, namely, consumption, a riskless asset, and a risky asset.  We assume that the individual invests in order to minimize the expectation of some nonincreasing, nonnegative function of her lifetime minimum wealth.  For further motivation of this problem, see Browne (1995, 1997, 1999a, b), Milevsky, Ho, and Robinson (1997), Hipp and Plum (2000), Hipp and Taksar (2000), Milevsky and Robinson (2000), Schmidli (2001), Young (2004), and Milevsky, Moore, and Young (2006).
The individual consumes at a Lipschitz continuous rate $c(w) \ge 0$, in which $w$ is her current wealth.  We assume that the individual invests in a riskless asset whose price at time $t$, $X_t$, follows the deterministic process $dX_t = rX_t dt, X_0 = x > 0$, for some fixed rate of interest $r > 0$. Also, the individual invests in a risky asset whose price at time $t$, $S_t$, follows geometric Brownian motion given by

$$dS_t = \mu S_t dt + \sigma S_t dB_t, \quad S_0 = S > 0, \eqno(2.1)$$
\noindent in which $\mu > r$, $\sigma > 0$, and $B$ is a standard Brownian motion with respect to a filtration $\{ {\cal F}_t \}$ of a probability space $(\Omega, {\cal F}, {\bf P})$.
Let $W_t$ be the wealth at time $t$ of the individual, and let $\pi_t$ be the amount that the decision maker invests in the risky asset at that time.  It follows that the amount invested in the riskless asset is $W_t - \pi_t$, and wealth follows the process

$$dW_t = [rW_t + (\mu - r) \pi_t  - c(W_t)] dt + \sigma \pi_t dB_t, \quad W_0 = w. \eqno(2.2)$$

\noindent A process associated with this wealth process is the {\it minimum} wealth process.  Let $M_t$ denote the minimum wealth of the individual during $[0, t]$; that is,

$$M_t = \min \left[ \inf_{0 \le s \le t} W_s, \tilde M_0 \right], \eqno(2.3)$$

\noindent in which we allow the individual to have a financial ``past'' by including $\tilde M_0$, the minimum wealth that the individual experienced before time 0.

By {\it lifetime minimum wealth}, we mean the minimum wealth between time 0 and the random time $\tau_d$ that the individual dies; $M_{\tau_d}$ denotes the lifetime minimum wealth.  We assume that $\tau_d$ is exponentially distributed with parameter $\lambda$ (that is, with expected time until death equal to $1/\lambda$); this parameter is also known as the {\it hazard rate} of the individual.  We assume that $\tau_d$ is independent of the $\sigma$-algebra generated by the Brownian motion $B$.

Denote the minimum expectation of a nonincreasing, nonnegative function $f$ of the lifetime minimum wealth by $V^f(w, m)$.  We minimize with respect to the set of admissible investment strategies $\cal A$.  A strategy $\pi$ is {\it admissible} if it is ${\{ {\cal F}_t \}}$-progressively measurable (in which ${\cal F}_t$ is the augmentation of $\sigma(W_s: 0 \le s \le t)$) and if it satisfies the integrability condition $\int_0^t \pi_s^2 \, ds < \infty$, almost surely, for all $t \ge 0$.  

It follows that one can express $V^f$ by

$$V^f(w, m) = \inf_{\pi \in {\cal A}} {\bf E}^{w,m} \left[f(M_{\tau_d}) \right], \eqno(2.4)$$

\noindent  in which ${\bf E}^{w, m}$ denotes conditional expectation given $W_0 = w$ and $M_0 = m$. Note that we do not discount the penalty function $f$ by a factor such as $e^{-\eta \tau_d}$ for some $\eta > 0$.  If we were to do so, then the resulting problem would be equivalent to adding $\eta$ to the hazard rate $\lambda$.

We restrict $f$ in (2.4) to be such that $V^f(w, m)$ is finite for all $w \ge m$.  We also require that an admissible strategy satisfy ${\bf E}^{w,m}[f(M_{\tau_d})] < \infty$.  Effectively, this is not a restriction because any strategy for which  ${\bf E}^{w,m}[f(M_{\tau_d})] = \infty$ will not be optimal.

For examples of relevant functions $f$, consider the following:  If $f(m) = {\bf 1}_{\{m \le b\}}$, then $f(M_{\tau_d})$ indicates whether the individual's wealth has reached $b \in {\bf R}$ during one's life, and $V^f$ is the minimum probability of lifetime ruin with ruin level $b$ (Young, 2004).  If $f(m) = \max(b-m, 0)$, then $V^f$ is the minimum expected lifetime shortfall relative to $b \in {\bf R}$ (Bayraktar and Young, 2005).

The function $f$ given by $f(m) = {\bf 1}_{\{m \le b\}}$ is special because we can write any nonincreasing, nonnegative function as the pointwise limit of an increasing sequence of functions of the form

$$f_n(m) = a_0 + \sum_{i = 1}^{d_n} a_i {\bf 1}_{\{m \le b_i\}}, \eqno(2.5)$$

\noindent in which $a_0 \ge 0$, $a_i > 0$ for $i = 1, 2, \dots, d_n$, and $b_1 > b_2 > \dots > b_{d_n}$; see, for example, Royden (1968).

Therefore, we begin by considering $V^f$ when $f(m) = {\bf 1}_{\{m \le b\}}$ for $b \in \bf R$.  In this case, we write $V(w, m; b)$ for $V^f(w, m)$, and

$$V(w, m; b) = \cases{1, &if $m \le b$; \cr
\psi(w; b), &if $m > b$;} \eqno(2.6)$$

\noindent in which $\psi(w; b)$ is the probability that wealth reaches $b$ before the individual dies, given that it has not done so to date.  In the next section, we provide a verification lemma for $\psi$.

For every $\alpha \in \bf R$, we associate a second-order differential operator ${\cal L}^{\alpha}$ with the minimization problem in (2.4) as follows:  For an open set $G \subset {\bf R}$ and for $v \in C^2(G)$, define the function ${\cal L}^{\alpha} v: G \rightarrow \bf R$ by

$${\cal L}^{\alpha} v(w) = \left[ rw + (\mu - r) \alpha - c(w) \right] v'(w) + {1 \over 2} \sigma^2 \alpha^2 v''(w) - \lambda v(w). \eqno(2.7)$$

\noindent We use ${\cal L}^{\alpha}$ in the next section to characterize $\psi$ in (2.6).

\subsect{2.2. Representation of the Value Function $V^f$}

Heinricher and Stockbridge (1991) provide a verification lemma for the $L^{\infty}$ control of a diffusion, whose drift and volatility depend on the control and satisfy some regularity assumptions: the volatility is bounded, and the derivatives of the volatility and drift with respect to their variables exist and are bounded.  (The boundedness of the value function is necessary for uniform integrability, which is why they require the boundedness of the volatility.)  Apart from this, they assume that the control problem ends as soon as the diffusion hits a prespecified level.  With these assumptions, they provide a verification lemma for the value function, which requires it to be smooth and satisfy a polynomial growth condition.  Note that the wealth function in (2.2) does not have a bounded volatility since the control space might be unbounded, and our problem  involves a diffusion with killing.  Specifically, the time of death of the individual kills the diffusion.

By using similar arguments to Heinricher and Stockbridge (1991), one can directly derive a verification lemma for $V^f$ in (2.4).  For such a verification lemma, we would need to require $V^f$ to be smooth and bounded, which turns out to be true if $f$ is smooth and bounded.  In this paper, however, the only assumptions we make on $f$ are that it is nonnegative and nonincreasing.  As a result, $V^f$ might not be smooth with respect to $m$; see (2.6) for an example in which $V^f$ is not differentiable with respect to $m$.  Instead of forcing a verification lemma for $V^f$, we derive one for $\psi$, as given in (2.6), via an argument that avoids the usual verification argument for $L^\infty$ control problems that exhibit a Neumann condition.  Then, we are able to show that one can represent $V^f$ in terms of $\psi$ by showing that the optimal strategy $\pi^*$ for (2.4) does not depend on $f$.  Therefore, our sufficient conditions for a function to be equal to the value function (2.4) are not as restrictive as Heinricher and Stockbridge (1991).

Barles, Daher, and Romano (1994), on the other hand, show that the value function of an $L^{\infty}$ control problem is a viscosity solution of the dynamic programming equation when the control resides on a compact space and the function corresponding to our $f$ is bounded.  As mentioned in the preceding paragraph, we represent $V^f$ explicitly in terms of $\psi$.  By using that expression, we observe that $V^f$ is a viscosity solution of the HJB equation with Neumann boundary condition without the compactness assumption of the control space and boundedness of $f$; see Corollary 2.7 below.

Define $w^s \equiv \inf \{ w: c(w) = rw, w \ge 0 \}$, with the convention that $\inf \phi = \infty$.  Note that when wealth equals $w^s$, then the individual can place all her wealth in the riskless asset and consume at the rate $c(w^s)$ for the remainder of her life without the risk of bankruptcy.  One can think of $w^s$ as the ``safe'' level.  Assume that the bankruptcy level $b$ is less than the safe level $w^s$.  Additionally, we have the unstated condition that the initial minimum wealth $M_0 = m > b$ whenever we discuss $\psi$. If $M_0 = m \le b$, then $\psi$ is undefined.

We have the following verification lemma for the probability of ruin $\psi$ whose proof is in the Appendix; for a general reference, see Fleming and Soner (1993).

\lem{2.1} {{\bf (Verification lemma)} Suppose $v$ is a nonincreasing, convex function from $[b, \infty)$ to $[0, 1]$ and suppose $\beta$ is a function from $[b, w^s)$ to ${\bf R}$ that satisfy the following conditions:
\item{$(i)$} $v$ is $C^2$ on $(b, \infty),$ except possibly at $w^s$ where it is $C^1;$
\item{$(ii)$} $\beta \in {\cal A},$ in which the investment strategy $\beta$ is defined by $\beta_t = \beta(W_t);$
\item{$(iii)$} ${\cal L}^{\alpha} v(w) \ge 0,$ for $\alpha \in {\bf R};$
\item{$(iv)$} ${\cal L}^{\beta(w)} v(w) = 0,$ for $w \in (b, w^s);$
\item{$(v)$} $v(b) = 1,$ and $v(w) = 0$ for $w \ge w^s$.}

{\it Under the above conditions, the minimum probability of the lifetime ruin $\psi$ is given by}
$$\psi(w; b) = v(w), \quad w \ge b, \eqno(2.8)$$

\noindent {\it and the optimal investment strategy in the risky asset $\pi^*$ is given by}
$$\pi^*(w; b) = \beta(w), \quad w \in (b, w^s). \eqno(2.9)$$

\medskip

From Lemma 2.1, we deduce that if we find a solution to the HJB equation embodied in that lemma, then that solution is the minimum probability of ruin $\psi$ as defined in (2.6).  Therefore, without abusing notation too greatly, we also represent the smooth, decreasing, convex solution on $[b, w^s)$ of the following HJB equation by $\psi$:

$$\left\{ \eqalign{&  \lambda \psi(w; b) = (rw - c(w)) \psi'(w; b) + \min_\pi \left[ (\mu - r) \pi \psi'(w; b) + {1 \over 2} \sigma^2 \pi^2 \psi''(w; b) \right], \cr
&  \psi(b; b) = 1, \quad \psi(w^s; b) = 0.} \right. \eqno(2.10)$$

\noindent If $w^s = \infty$, then $\psi(w^s; b) = 0$ means that $\lim_{w \rightarrow \infty} \psi(w; b) = 0$.  One can show via standard techniques (Zariphopoulou, 1994) that $\psi$ given in (2.6) is a viscosity solution of (2.10).  Lemma 2.1 then tells us that if we find a smooth solution $v$ of (2.10), then $\psi = v$ is a smooth solution of (2.10).  For the sake of brevity, we do not determine the general conditions under which (2.10) has a smooth solution, although we construct such a solution in Section 4.1 for a case pertinent to this paper.

A key observation is given in the following proposition:

\prop{2.2} {Suppose the following problem has a convex solution $h \in C^1({\bf R}) \cap C^2({\bf R} - \{w^s\}):$
$$\left\{ \eqalign{&  \lambda h(w) = (rw - c(w)) h'(w) + \min_\pi \left[ (\mu - r) \pi h'(w) + {1 \over 2} \sigma^2 \pi^2 h''(w) \right], \quad -\infty < w < w^s, \cr
& h(0) = 1, \hbox{and } h(w) = 0 \hbox{ for } w \ge w^s.} \right. \eqno(2.11)$$
Then,
$$\psi(w; b) = h(w)/h(b), \quad w \in [b, \infty), \eqno(2.12)$$
and
$$\pi^*(w; b) = -{\mu - r \over \sigma^2} {h'(w) \over h''(w)}, \quad w \in [b, w^s). \eqno(2.13)$$
Additionally, the solution of $(2.11)$ is unique on $\bf R$.}

\medskip

\pf  Define $v$ by $v(w) = h(w)/h(b)$ on $[b, \infty)$, and define $\beta(w) = -{\mu - r \over \sigma^2} {h'(w) \over h''(w)}$ on $[b, w^s)$.   It is clear that $v$ and $\beta$ satisfy the conditions of Lemma 2.1, so (2.12) and (2.13) follow.  For any $b < 0$, note that $h(w) = \psi(w; b)/\psi(0; b)$ solves (2.11) on $[b, \infty)$.  By Lemma 2.1, $\psi$ is unique on $[b, \infty)$; therefore, $h$ is unique on $[b, \infty)$.  Because $b < 0$ is arbitrary, the solution of (2.11) is unique on $\bf R$.   $\square$

\medskip

The boundary condition $h(0) = 1$ is a normalizing condition.  One could set $h(0) = y$ for any $y > 0$ and still obtain that $\psi(w; b) = h(w)/h(b)$.  Note that for our choice of $h(0) = 1$, we have $h(w) = \psi(w; 0)$ for $w \ge 0$.

Now that we have a simpler representation of $\psi$ given in (2.12) via the solution $h$ of (2.11), we focus henceforth on the function $h$.  We proceed to obtain a representation of $V^f$ in terms of $h$; see Theorem 2.5 below.  First, we have an immediate corollary that follows from the expression for $\pi^*$ in (2.13).

\cor{2.3} {Suppose $(2.11)$ has a convex solution $h \in C^1({\bf R}) \cap C^2({\bf R} - \{w^s\})$; then, the optimal investment strategy $\pi^*$ is independent of the ruin level $b$.}

\medskip

In the remainder of this section, we assume that a convex solution $h \in C^1({\bf R}) \cap C^2({\bf R} - \{w^s\})$ to (2.11) exists.  In Section 4.1, we find such a solution for a special case of the rate of consumption $c$.  Next, suppose that $f_n$ is given by (2.5).  Then, we have the following lemma used in the proof of Theorem 2.5 below.  It follows easily from Proposition 2.2 and Corollary 2.3, so we omit its proof.

\lem{2.4} {For $f_n$ in $(2.5)$, we have that $V^{f_n}$ is given by}

$$V^{f_n} (w, m) = a_0 + \cases{\sum_{i = 1}^{d_n} a_i, & if $m \le b_{d_n}$, \cr
\sum_{i = 1}^{d_n-1} a_i + a_{d_n} h(w)/h(b_{d_n}), & if $b_{d_n} < m \le b_{d_n -1}$, \cr
\sum_{i = 1}^{d_n-2} a_i + \sum_{i = d_n-1}^{d_n} a_i h(w)/h(b_i), & if $b_{d_n-1} < m \le b_{d_n-2}$, \cr
\dots \cr
a_1 + \sum_{i = 2}^{d_n} a_i h(w)/h(b_i), & if $b_2 < m \le b_1$, \cr
\sum_{i = 1}^{d_n} a_i h(w)/h(b_i), & if $m > b_1$.}  \eqno(2.14)$$

\noindent {\it The corresponding optimal investment strategy $\pi^*$ is given by $(2.13)$.}

\medskip

The following theorem is one of our main contributions.  It tells us that the optimal investment strategy $\pi^f$ corresponding to $V^f$ for a nonincreasing, nonnegative function $f$ is $\pi^*$ in (2.13).  In particular, the optimal investment strategy corresponding to $V^f$ is independent of $f$.

\th{2.5} {Let $f$ be a nonincreasing, nonnegative function for which $V^f(w, m)$ is finite for all $w \in [m, w^s)$.  Then, $V^f(w, m) = {\bf E}^{w, m} \left[ f(M^{\pi^*}_{\tau_d}) \right],$ in which the optimal investment strategy $\pi^*$ is given by $(2.13)$.}

\medskip

\pf Let $\{ f_n \}$ be an increasing sequence of functions of the form in (2.5) such that $f(m) = \lim_n f_n(m)$ for each $m \in {\bf R}$.  Thus,

$$\eqalign{{\bf E}^{w, m} \left[ f(M^{\pi^*}_{\tau_d}) \right] &= {\bf E}^{w, m} \left[ \lim_n f_n(M^{\pi^*}_{\tau_d}) \right] = \lim_n {\bf E}^{w, m} \left[ f_n(M^{\pi^*}_{\tau_d}) \right] \cr
&= \lim_n \inf_{\pi \in {\cal A}} {\bf E}^{w, m} \left[ f_n(M^\pi_{\tau_d}) \right] \ \le \lim_n {\bf E}^{w, m} \left[ f_n(M^\pi_{\tau_d}) \right] \cr
&\le {\bf E}^{w, m} \left[ f(M^\pi_{\tau_d}) \right],} \eqno(2.15)$$

\noindent {in which the second equality follows from the Monotone Convergence Theorem (Royden, 1968, Theorem, 4.9), and in which the last inequality follows from the fact that $f_n \le f$.  Now, by taking the infinimum over admissible strategies $\pi$, we have that ${\bf E}^{w, m} \left[ f(M^{\pi^*}_{\tau_d}) \right] \le \inf_{\pi \in {\cal A}} {\bf E}^{w, m} \left[ f(M^\pi_{\tau_d}) \right]$.  Thus, the theorem follows.}  $\square$

\medskip

From Theorem 2.5, we obtain a useful representation of $V^f$ in terms of the probability of ruin $\psi$, or equivalently in terms of $h$ given by (2.11).  Recall that we assume that a convex solution $h \in C^1({\bf R}) \cap C^2({\bf R} - \{w^s\})$ of (2.11) exists.

\cor{2.6} {We can express $V^f$ for $m \le w < w^s$ as follows:
$$\eqalign{V^f(w, m) & = f(m) \left( 1 - {h(w) \over h(m)} \right) + \int_{-\infty}^m f(x) {\partial \over \partial x} {h(w) \over h(x)} dx \cr
&= f(m) - \int_{-\infty}^m f'(x) {h(w) \over h(x)} dx.}  \eqno(2.16)$$}

\medskip

\pf  Because $\pi^*$ in (2.13) is the optimal investment strategy for $V^f$, it follows that

$$V^f(w, m) = {\bf E}^{w, m} [f(M^*_{\tau_d})] = \int_{-\infty}^\infty f(x) d{\bf P}^{w, m} [M^*_{\tau_d} \le x], \eqno(2.17)$$

\noindent in which $M^*$ denotes the optimally-controlled minimum wealth $M^{\pi^*}$.  Note that ${\bf P}^{w, m} [M^*_{\tau_d} \le x] = V(w, m; x)$ from (2.6), in which ${\bf P}^{w, m}$ denotes the conditional probability given $W_0 = w$ and $M_0 = m$.  Thus,

$$\int_{-\infty}^\infty f(x) d{\bf P}^{w, m} [M^*_{\tau_d} \le x] = f(m) \left( 1 - \psi(w; m) \right) + \int_{-\infty}^m f(x) {\partial \over \partial x} \psi(w; x) dx, \eqno(2.18)$$

\noindent which is the first line of (2.16).  The second line of (2.16) follows from integration by parts if $f$ is differentiable.  If $f$ is not differentiable, then we interpret the derivative $f'$ in the sense of a distribution (Al-Gwaiz, 1992).  $\square$

\medskip

From (2.16), we deduce that $V^f$ is a viscosity solution of an associated HJB equation.  We state this formally and without proof in the following corollary, which is similar to a result of Barles, Daher, and Romano (1994).

\cor{2.7} {$V^f$ in $(2.4)$ is a viscosity solution of}
$$\left\{ \eqalign{&  \lambda V^f(w, m) = (rw - c(w)) V^f_w(w, m) + \min_\pi \left[ (\mu - r) \pi V^f_w(w, m) + {1 \over 2} \sigma^2 \pi^2 V^f_{ww}(w, m) \right], \cr
&  V^f(w^s, m) = f(m), \quad V^f_m(m, m) = 0.} \right. \eqno(2.19)$$

\medskip

We next turn our attention to finding when the optimal investment strategy for the problem considered in this section equals the optimal investment strategy for the problem of maximizing utility of lifetime consumption.

\sect{3. Relation with Maximizing Utility of Consumption in a Black-Scholes Market}

In this section, we explore when the optimal investment strategy for the model in Section 2 coincides with the one for a model in which an individual seeks to maximize expected discounted utility of lifetime consumption.  We assume that the optimal consumption strategy arising from maximizing utility of consumption is the consumption that the individual then follows in the minimum wealth problem.

In Section 3.1, we state the problem of maximizing the expected discounted utility of lifetime consumption, and we present the corresponding HJB equation that the value function solves.  These results are well-documented in the literature, so we cover them only briefly; see, for example, Merton (1992) and Karatzas and Shreve (1998, Chapter 3).  In Section 3.2, we equate the optimal investment strategy when minimizing a function of lifetime minimum wealth, $\pi^*$ in (2.13), with the optimal investment strategy when maximizing utility of lifetime consumption, $\pi^u$.  We determine the resulting relationship between the corresponding value functions and what consumption function is implied for the two problems by assuming that the rate of consumption is the same for both.  We show that the investment strategies equal only when the consumption function is linear.

\subsect{3.1. Maximizing Expected Utility of Lifetime Consumption}

We assume that the individual invests in a financial market, as described in Section 2.1, and that the individual consumes in order to maximize her expected discounted utility of consumption during her life.  Similar to Section 2.1, wealth follows the process given in (2.2), except that we now take the rate of consumption as a control; it is not given exogenously as in (2.2).

We assume that the individual maximizes her utility of consumption over admissible investment and consumption strategies ${\cal A}'$, in which an investment strategy $\pi$ is admissible if it is $\{ {\cal F}_t \}$-progressively measurable and if it satisfies the integrability condition $\int_0^t \pi_s^2 \, ds < \infty$, almost surely, for all $t$.  A consumption strategy $c$ is admissible if it is nonnegative, if it is $\{ {\cal F}_t \}$-progressively measurable, and if it satisfies the integrability condition $\int_0^t c_s \, ds < \infty$, almost surely, for all $t$.

Let $u$ represent the individual's utility of consumption.  As in Karatzas and Shreve (1998, Section 3.4), we assume that $u: {\bf R} \rightarrow [-\infty, \infty)$ is a concave, nondecreasing, upper semicontinuous function that satisfies: (1) The half-line $dom(u) \equiv \{x \in {\bf R}: u(c) > -\infty \}$ is a nonempty subset of $[0, \infty)$; and (2) $u'$ is continuous, positive, strictly decreasing on the interior of $dom(u)$, and $\lim_{x \rightarrow \infty} u'(x) = 0$.

We define the the value function $V^u$ by

$$V^u(w) = \sup_{(\pi, c) \in {\cal A}'} {\bf E}^w \left[ \int_0^{\tau_d} e^{-\rho s} u(c_s) ds \right], \eqno(3.1)$$

\noindent in which $\rho \ge 0$ is the individual subjective discount rate, a measure of impatience.  We will see in Theorem 3.1 below that the individual subjective discount rate $\rho$ is the slope of the linear consumption rate that results when we equate the investment strategies from maximizing utility and minimizing a function of minimum wealth.   As in Section 2, we assume that $\tau_d$ is distributed exponentially with mean $1/\lambda$ and that $\tau_d$ is independent of $\sigma$-algebra generated by the Brownian motion $B$.  Because $u$ is concave and because wealth is linear in the controls, $V^u$ is concave.

It follows from Karatzas and Shreve  (1998, Theorem 3.9.20) that under a finiteness condition on the model (see their equation (3.9.25)), $V^u$ exists and solves the following HJB equation for $w > w_\infty \equiv \inf\{ x \in {\bf R}: V^u(x) > -\infty \}$:

$$(\rho + \lambda) V^u = rw \cdot (V^u)' + \max_{c \ge 0} \left[ u(c) - c \cdot (V^u)' \right] + \max_\pi \left[ (\mu - r) \pi \cdot (V^u)' + {1 \over 2} \sigma^2 \pi^2 \cdot (V^u)'' \right], \eqno(3.2)$$

\noindent with corresponding optimal investment strategy $\pi^u$ given in feedback form by

$$\pi^u(w) = -{\mu - r \over \sigma} {(V^u)'(w) \over (V^u)''(w)},  \eqno(3.3)$$

\noindent and with the optimal consumption strategy $c^u$ solving

$$u'(c^u(w)) =  (V^u)'(w). \eqno(3.4)$$

\subsect{3.2. Equating the Optimal Investment Strategies}

In this section, we equate the optimal investment strategies for the problems of minimizing lifetime ruin and maximizing utility of consumption and determine that the consumption strategy is a linear function of wealth, if the minimizer of lifetime ruin follows the consumption dictated by the maximizer of utility.  We make this more precise in the following theorem:

\th{3.1} {Suppose that $(2.11)$ has a convex solution $h \in C^1({\bf R}) \cap C^2({\bf R} - \{w^s\})$, with corresponding optimal investment strategy $\pi^*$ given in $(2.13)$. Suppose that $\pi^*(w) = \pi^u(w)$ on $(\kappa, w^s)$, where $\kappa \in (\max(0, w_\infty), w^s)$, in which $w_\infty = \inf\{ x \in {\bf R}: V^u(x) > -\infty \}$ and $w^s = \inf \{ w: c(w) = rw, w \ge 0 \}$ is the safe level.  Moreover, suppose the rate consumption $c$ in the wealth process $(2.2)$ is given by $(3.4)$.  Then, the rate of consumption is linear on $(\kappa, w^s)$ with slope equal to the personal discount rate $\rho$.}

\medskip

\pf  When we equate $\pi^*$ from (2.13) with $\pi^u$ from (3.3), we obtain

$${h'(w) \over h''(w)} = {(V^u)'(w) \over (V^u)''(w)}, \eqno(3.5)$$

\noindent from which it follows that

$$h'(w) = k \cdot (V^u)'(w), \eqno(3.6)$$

\noindent for some constant $k < 0$.

Thus,

$$\left[ (rw - c(w)) - \delta {(V^u)'(w) \over (V^u)''(w)} \right] k \cdot (V^u)'(w)= \left[ (rw - c(w)) - \delta {h'(w) \over h''(w)} \right] h'(w).  \eqno(3.7)$$

\noindent Rewrite the left-hand side of (3.7) by using equation (3.2), and rewrite the right-hand side by using equation (2.11) to obtain

$$(\rho + \lambda) \, k \, V^u(w) - k \, u(c(w)) = \lambda h(w).  \eqno(3.8)$$

Note that we can solve (3.3) for $V^u$ in terms of $\pi^u$.  Indeed, for $w \in (\kappa, w^s)$,

$$V^u(w) = V^u(w_0) + (V^u)'(w_0) \int_{w_0}^w \exp \left\{ -{\mu - r \over \sigma^2} \int_{w_0}^v {dz \over \pi^u(z)} \right\} dv,  \eqno(3.9)$$

\noindent in which $w_0$ is some point in $(\kappa, w^s)$.  We have a similar expression for $h$ in terms of $\pi^*$.  For the remainder of the proof, we write $\pi$ for both $\pi^u$ and $\pi^*$ because we assume $\pi^u = \pi^*$ on $(\kappa, w^s)$.

Note that $c$ in (3.4) has a (strictly) positive derivative on $(\kappa, w^s)$; indeed, $c'(w) = (V^u)''(w)/u''(c(w)) > 0$ because both $V^u$ and $u$ are strictly concave on $(\kappa, w^s)$.  Therefore, we can solve (3.4) for $u$ in terms of $\pi$ and obtain that for $w \in (\kappa, w^s)$,

$$u(c(w)) = u(c(w_0)) + (V^u)'(w_0) \int_{w_0}^w \exp \left\{ -{\mu - r \over \sigma^2} \int_{w_0}^v {dz \over \pi(z)} \right\} \cdot c'(v) dv.  \eqno(3.10)$$

\noindent Substitute expressions (3.9) and (3.10) into (3.8) to get

$$\eqalign{& (\rho + \lambda) k \left[ V^u(w_0) + (V^u)'(w_0) \int_{w_0}^w \exp \left\{ -{\mu - r \over \sigma^2} \int_{w_0}^v {dz \over \pi(z)} \right\} dv \right] \cr
& \qquad - k \left[ u(c(w_0)) + (V^u)'(w_0) \int_{w_0}^w \exp \left\{ -{\mu - r \over \sigma^2} \int_{w_0}^v {dz \over \pi(z)} \right\} \cdot c'(v) dv \right] \cr
& =  \lambda \left[ h(w_0) + h'(w_0) \int_{w_0}^w \exp \left\{ -{\mu - r \over \sigma^2} \int_{w_0}^v {dz \over \pi(z)} \right\} dv \right].} \eqno(3.11)$$

\noindent Because (3.8) holds for $w = w_0$, the ``constant'' terms in (3.11) cancel, and we are left with

$$\int_{w_0}^w \exp \left\{ -{\mu - r \over \sigma^2} \int_{w_0}^v {dz \over \pi(z)} \right\} \cdot (\rho - c'(v)) dv = 0,  \eqno(3.12)$$

\noindent from which it follows that for $w \in (\kappa, w^s)$,

$$c'(w) = \rho;  \eqno(3.13)$$

\noindent that is, $c$ is a linear function of wealth on $(\kappa, w^s)$.  We have, thus, proved the theorem.  $\square$

\medskip

Note that if we were to replace $f(M_{\tau_d})$ with $e^{-\eta \tau_d} f(M_{\tau_d})$ in the definition of the value function $V^f$ in (2.4), then (3.13) would change to $c'(w) = \rho - \eta$, and the analysis that follows in this paper would go through with $\rho$ replaced by $\rho - \eta$.  For $c'(w) > 0$, we would require $\eta < \rho$.  This assumption makes sense because $\rho$ applies to the utility of current consumption, whereas $\eta$ in $e^{-\eta \tau_d}$ discounts from the time of death, but the minimum wealth almost surely attains $M_{\tau_d}$ before time $\tau_d$.

Theorem 3.1 does not hold for more general financial models.  For example, if $r$ and $\rho$ are deterministic functions of time, then we can no longer conclude that by equating the investment strategies that the consumption function is linear with respect to wealth with slope equal to the discount rate $\rho$.  The reason is that $k$ in equation (3.6) is no longer a constant.  It will be a function of time, so the proof does not go through as before.

In Section 4, we obtain a type of converse of Theorem 3.1.  We first show that if (3.13) holds on $(\kappa, w^s)$, then the optimal investment strategy for minimizing the probability of ruin is a linear function of wealth.  Then, we show that the utility function implied by these consumption and investment strategies exhibits hyperbolic risk aversion on $(\kappa, w^s)$.

\sect{4. Piecewise Linear Consumption}

Theorem 3.1 states that the {\it only} time the investment strategies for minimizing a function of lifetime minimum wealth and for maximizing utility of consumption coincide is when the rate of consumption is a linear function of  wealth (also assuming that the rates of consumption coincide).  Therefore, in this section, we assume that $c'(w) = \rho$ on $(\kappa, w^s)$, as in (3.13).  In order to keep the problem of minimizing a function of lifetime minimum wealth from being trivial, we suppose that $c(w)$ is bounded below by a positive constant for all $m \le w < w^s$.  We use a simple case of such a bound by assuming that $c$ equals a positive constant for $w \le \kappa$.  Specifically, let the rate of consumption be given by

$$c(w) = (\bar c + \rho \kappa) + \rho(w - \kappa)_+, \eqno(4.1)$$

\noindent with $\rho > 0$, $\kappa > 0$, and $\bar c$ such that $\bar c + \rho \kappa > 0$.  In fact, we assume that $\bar c + \rho \kappa > r \kappa$.  If we were to allow $\bar c + \rho \kappa \le r \kappa$, then $w^s \le \kappa$.  In this case, the rate of consumption would equal the constant $\bar c + \rho \kappa$ on $(-\infty, w^s)$.  Young (2004) solves this problem for minimizing the probability of lifetime ruin, and Bayraktar and Young (2005) solve it for minimizing expected lifetime shortfall.

In Section 4.1, we calculate the minimum probability of lifetime ruin under the consumption function given by (4.1).  Additionally, we learn that in general, if $c'(w) = \rho$ on $(\kappa, w^s)$, then the corresponding optimal investment strategy is linear on $(\kappa, w^s)$.  In Section 4.2, we set the optimal investment strategy $\pi^u$ when maximizing utility equal to this linear function on $(\kappa, w^s)$.  We show that this, in turn, implies that the utility function given by (3.10) exhibits hyperbolic absolute risk aversion (HARA) on $(\kappa, w^s)$.  Therefore, the only investors that behave similarly under the two problems of minimizing a function of lifetime minimum wealth or maximizing utility of consumption are those for which the utility is HARA.

\subsect{4.1.  Minimizing the Probability of Lifetime Ruin}

In Proposition 2.2, we showed that the minimum probability of lifetime ruin $\psi(w; b)$ with ruin level $b$ equals $h(w)/h(b)$, in which $h$ solves (2.11).  We proceed by solving this nonlinear equation for $h$.  We assume {\it a priori} that the solution is nonincreasing, convex, and $C^2$ on $\bf R$ (except possibly at $w^s$ where it is $C^1$) and by construction show that our solution $h$ satisfies these properties.   Because we are looking for a convex solution, we consider the Legendre transform $\tilde h$ of $h$ for $w < w^s$ defined by

$$\tilde h(y) \equiv \min_{w < w^s} \left[ h(w) + w y \right];  \eqno(4.2)$$

\noindent see Karatzas and Shreve (1998, Chapter 3).  The function $\tilde h$ is also called the concave dual of $h$.  Note that we can recover $h$ from $\tilde h$ by

$$h(w) = \max_{y > 0} [\tilde h(y) - wy]. \eqno(4.3)$$

\noindent The minimizing value of $w$ in (4.2) equals $I(-y) = \tilde h'(y)$, in which $I$ is the inverse function of $h'$.  Therefore, the maximizing value of $y$ in (4.3) equals $-h'(w)$.

Substitute $w = I(-y) = \tilde h'(y)$ in equation (2.11) to obtain

$$\lambda \tilde h(y) + (r - \lambda) y \tilde h'(y) - \delta y^2 \tilde h''(y) = c(\tilde h'(y)) y, \eqno(4.4)$$

\noindent in which $\delta = {1 \over 2} \left( {\mu - r \over \sigma} \right)^2$.  Note that (4.4) is a linear differential equation if and only if $c(w)$ is a linear function of $w$.

Recall that $h \in C^1({\bf R}) \cap C^2({\bf R} - \{w^s\})$; in particular, $h$ is $C^1$ at $w = w^s$, so $h'(w_s) = 0$ if $w^s < \infty$.  If $w^s = \infty$, then we also make use of the assumption that $h$ is decreasing, nonnegative, and convex on $\bf R$ to assert that $h'(w^s) = 0$.  It follows that the values of $y$ that correspond to $w = 0$ and $w = w^s$ are $y_0 = -h'(0)$ and $y_s = -h'(w^s) = 0$, respectively.  In terms of $\tilde h$, we can write these expressions as

$$\tilde h'(0) = w^s, \hbox{ and } \tilde h'(y_0) = 0, \eqno(4.5)$$

\noindent and the conditions $h(w^s) = 0$ and $h(0) = 1$ become

$$\tilde h(0) = 0, \hbox{ and } \tilde h(y_0) = 1. \eqno(4.6)$$

To solve (4.4), we begin by splitting its domain into two pieces:  One for $w < \kappa$, and the other for $\kappa < w < w^s$.  If $w < \kappa$, then $c(w) = \bar c + \rho \kappa$, and equation (4.4) becomes

$$\lambda \tilde h(y) + (r - \lambda) y \tilde h'(y) - \delta y^2 \tilde h''(y) = (\bar c + \rho \kappa) y, \eqno(4.7)$$

\noindent for $y > y_\kappa$, in which $y_\kappa = -h'(\kappa-)$, or equivalently

$$\tilde h'(y_\kappa+) = \kappa.  \eqno(4.8)$$

The general solution of (4.7) is

$$\tilde h(y) = D_1 y^{B_1} + D_2 y^{B_2} + {\bar c + \rho \kappa \over r} y,  \eqno(4.9)$$

\noindent in which $D_1$ and $D_2$ are constants to be determined by the boundary conditions, and $B_1 > 1$ and $B_2 < 0$ are the positive and negative solutions, respectively, of the following quadratic equation:

$$\delta B^2 - (r - \lambda + \delta) B - \lambda = 0. \eqno(4.10)$$

If $\kappa < w < w^s$, then $c(w) = \bar c + \rho w$, and equation (4.4) becomes

$$\lambda \tilde h(y) + (r - \rho - \lambda) y \tilde h'(y) - \delta y^2 \tilde h''(y) = \bar c y,  \eqno(4.11)$$

\noindent for $0 < y < y_\kappa$, in which $y_\kappa = -h'(\kappa+)$, or equivalently

$$\tilde h'(y_\kappa-) = \kappa.  \eqno(4.12)$$

\noindent That is, we require that $\tilde h'$ be continuous at $y = y_\kappa$.  Below, we also use that $\tilde h$ is continuous at $y = y_\kappa$ to determine $\tilde h$ for $y > 0$.

The solution of (4.11) depends on whether or not $\rho$ equals $r$.  Also, depending on whether or not $\rho$ is less than $r$, the boundary $w^s$ is finite or infinite, respectively.  At this point, we discuss each case ($\rho > r$, $\rho = r$, and $\rho < r$) individually.

\medskip

\noindent{\bf A.} $\rho > r$.  In this case, $w^s = \infty$, and the general solution of (4.11) is given by

$$\tilde h(y) = \hat D_1 y^{\hat B_1} + \hat D_2 y^{\hat B_2} - {\bar c \over \rho - r}y, \eqno(4.13)$$

\noindent in which $\hat D_1$ and $\hat D_2$ are constants to be determined by the boundary conditions, and $\hat B_1 \in (0, 1)$ and $\hat B_2 < 0$ are the positive and negative solutions, respectively, of the following quadratic equation:

$$\delta \hat B^2 - (r - \rho - \lambda + \delta)\hat B - \lambda = 0. \eqno(4.14)$$

The first part of equation (4.6), namely $\tilde h(0) = 0$, implies that $\hat D_2 = 0$.  Thus, it follows from (4.3) that on the interval $(\kappa, \infty)$, the function $h$ is given by

$$h(w) = K_1 \left(w + {\bar c \over \rho - r} \right)^d, \eqno(4.15)$$

\noindent for some constant $K_1$ that depends on $\hat D_1$ and $\hat B_1$ and for $d = \hat B_1/(\hat B_1 - 1) < 0$.  Specifically,

$$d = {1 \over 2(r - \rho)} \left[ (r - \rho + \lambda + \delta) + \sqrt{(r - \rho + \lambda + \delta)^2 + 4(\rho - r) \lambda} \right].  \eqno(4.16)$$

Because $h$ is convex with respect to $w$ on $(\kappa, \infty)$, we can obtain the optimal investment strategy from the first-order necessary condition in (2.11).  Thus, the optimal investment strategy on $(\kappa, \infty)$ equals

$$\pi^*(w) = {\mu - r \over \sigma^2}  {w + {\bar c \over \rho - r} \over 1 - d},  \eqno(4.17)$$

\noindent a linear function in wealth that increases with respect to wealth as $\rho > r$.

Next, we show that we can determine the unknown constants $y_0$, $y_\kappa$, $D_1$, $D_2$, and $\hat D_1$ from the boundary conditions at $y = y_0$ and $y = y_\kappa$.  For reference, we rewrite the conditions in terms of the unknown constants.  The second parts of equations (4.5) and (4.6), respectively, imply that

$$D_1 B_1 y_0^{B_1 - 1} + D_2 B_2 y_0^{B_2 - 1} + {\bar c + \rho \kappa \over r} = 0, \eqno(4.18)$$

\noindent and

$$D_1 y_0^{B_1} + D_2 y_0^{B_2} + {\bar c + \rho \kappa \over r} y_0 = 1.  \eqno(4.19)$$

Equations (4.8) and (4.12), respectively, imply that

$$D_1 B_1 y_\kappa^{B_1 - 1} + D_2 B_2 y_\kappa^{B_2 - 1} + {\bar c + \rho \kappa \over r} = \kappa, \eqno(4.20)$$

\noindent and

$$\hat D_1 \hat B_1 y_\kappa^{\hat B_1 - 1} - {\bar c  \over \rho - r} = \kappa. \eqno(4.21)$$

\noindent Finally, continuity of $\tilde h$ at $y = y_\kappa$ implies that

$$\hat D_1 y_\kappa^{\hat B_1} - {\bar c \over \rho - r} y_\kappa = D_1 y_\kappa^{B_1} + D_2 y_\kappa^{B_2} + {\bar c + \rho \kappa \over r} y_\kappa. \eqno(4.22)$$

First, eliminate $\hat D_1$ from (4.21) and (4.22) to get an equation involving only $D_1$, $D_2$, and $y_\kappa$.  Solve the resulting equation and (4.20) for $D_1$ and $D_2$ in terms of $y_\kappa$.

$$D_1 = -{y_\kappa^{1 - B_1} \over \hat B_1 (B_1 - B_2)} \left[ {\bar c + \rho \kappa \over r} \hat B_1(1 - B_2) - \kappa (\hat B_1 - B_2) + {\bar c \over \rho - r} B_2 (1 - \hat B_1) \right] < 0, \eqno(4.23)$$

\noindent and

$$D_2 = -{y_\kappa^{1 - B_2} \over \hat B_1 (B_1 - B_2)} \left[ {\bar c + \rho \kappa \over r} \hat B_1(B_1 - 1) - \kappa (B_1 - \hat B_1) - {\bar c \over \rho - r} B_1 (1 - \hat B_1) \right]. \eqno(4.24)$$

Next, substitute (4.23) and (4.24) into (4.18) to get an equation for $y_0/y_\kappa$.

$$\eqalign{& {B_1 (y_0/y_\kappa)^{B_1 - 1} \over \hat B_1 (B_1 - B_2)} \left[ {\bar c + \rho \kappa \over r} \hat B_1(1 - B_2) - \kappa (\hat B_1 - B_2) + {\bar c \over \rho - r} B_2 (1 - \hat B_1) \right] \cr
& + {B_2 (y_0/y_\kappa)^{B_2 - 1} \over \hat B_1 (B_1 - B_2)} \left[ {\bar c + \rho \kappa \over r} \hat B_1(B_1 - 1) - \kappa (B_1 - \hat B_1) - {\bar c \over \rho - r} B_1 (1 - \hat B_1) \right] = {\bar c + \rho \kappa \over r}.} \eqno(4.25)$$

\noindent That (4.25) has a unique solution $y_0/y_\kappa > 1$ follows from the following three observations:

\item{(1) } If we set $y_0/y_\kappa = 1$, the left-hand side of (4.25) equals $(\bar c + \rho \kappa)/r - \kappa < (\bar c + \rho \kappa)/r$.

\item{(2) } As $y_0/y_\kappa$ approaches $\infty$, the left-hand side of (4.25) also approaches $\infty$ because $D_1 < 0$.

\item{(3) } The left-hand side of (4.25) is strictly increasing with respect to $y_0/y_\kappa$ for $y_0/y_\kappa > 1$.

\smallskip

Next, substitute (4.23) and (4.24) into (4.19) to get an equation for $y_0$ in terms of $y_0/y_\kappa$. Specifically,

$$\eqalign{& -{(y_0/y_\kappa)^{B_1 - 1} \over \hat B_1 (B_1 - B_2)} \left[ {\bar c + \rho \kappa \over r} \hat B_1(1 - B_2) - \kappa (\hat B_1 - B_2) + {\bar c \over \rho - r} B_2 (1 - \hat B_1) \right] + {\bar c + \rho \kappa \over r}  \cr
& - {(y_0/y_\kappa)^{B_2 - 1} \over \hat B_1 (B_1 - B_2)} \left[ {\bar c + \rho \kappa \over r} \hat B_1(B_1 - 1) - \kappa (B_1 - \hat B_1) - {\bar c \over \rho - r} B_1 (1 - \hat B_1) \right] = {1 \over y_0}.} \eqno(4.26)$$

\noindent Solve (4.26) for $y_0$; then, $y_\kappa$ is given by ${y_0 \over y_0/y_\kappa}$.  Next, compute the expressions in (4.23) and (4.24) to determine $D_1$ and $D_2$, respectively.  Finally, get $\hat D_1$ from either (4.21) or (4.22).

Once we have $\tilde h$, we can recover $h$ from (4.3).  Indeed, for $w < \kappa$, or equivalently $y > y_\kappa$,

$$h(w) = D_1 y^{B_1} + D_2 y^{B_2} + \left( {\bar c + \rho \kappa \over r} - w \right)y, \eqno(4.27)$$

\noindent in which $y$ solves

$$D_1 B_1 y^{B_1 - 1} + D_2 B_2 y^{B_2 - 1} + {\bar c + \rho \kappa \over r} - w = 0.  \eqno(4.28)$$

\noindent Given $w < \kappa$, solve equation (4.28) for $y$ and substitute into (4.27) to get $h(w)$.  Through a fair amount of algebra, one can show that the left-hand side of (4.28) decreases with respect to $y$ for $y > y_\kappa$; showing this is equivalent to showing that the left-hand side of (4.25) increases with respect to $y_0/y_\kappa > 1$. It follows that (4.28) has a unique solution for any $w < \kappa$.  Incidentally, this observation confirms that $\tilde h$ is concave and that, thereby, $h$ is convex on $(-\infty, \kappa)$.  Thus, $h \in C^2({\bf R})$ is nonincreasing and convex.

Earlier we observed that the optimal investment strategy $\pi^*$ is given by (4.17) for $w > \kappa$.  Because $h$ is convex with respect to $w$, we can obtain the optimal investment strategy from the first-order necessary condition in (2.11).  Thus, for $w < \kappa$,

$$\pi^*(w) = -{\mu - r \over \sigma^2} \left[ D_1 B_1(B_1 - 1) y^{B_1 - 1} + D_2 B_2(B_2 - 1) y^{B_2 - 1} \right], \eqno(4.29)$$

\noindent in which $y$ solves equation (4.28). Via tedious algebra, one can show that $\pi^*$ in (4.29) {\it decreases} with respect to wealth, although once $w > \kappa$, we know from (4.17) that $\pi^*$ increases with wealth.  We find this myopia on the part of the investor rather interesting.  Another surprising result is that $\pi^*$ decreases as $\rho$ increases for $w < \kappa$, although not surprisingly, once wealth is large enough, $\pi^*$ increases as $\rho$ increases.

\medskip

\noindent{\bf B.} $\rho = r$.  In this case, $w^s = \infty$ and the general solution of (4.11) is given by

$$\tilde h(y) = \hat D_1 y + \hat D_2 y^{-{\lambda \over \delta}} - {\bar c \over \delta + \lambda}y \ln y, \eqno(4.30)$$

\noindent in which $\hat D_1$ and $\hat D_2$ are constants to be determined by the boundary conditions.  As in Case A, $\hat D_2 = 0$, so we can determine that $h$ for $\kappa < w < \infty$ is given by

$$h(w) = K_2 \exp \left( - {\delta + \lambda \over \bar c} w \right),  \eqno(4.31)$$

\noindent for some constant $K_2$ that depends on $\hat D_1$.  From (4.31), it follows that the optimal investment strategy for $\kappa < w < \infty$ is given by

$$\pi^*(w) = {\mu - r \over \sigma^2} {\bar c \over \delta + \lambda}, \eqno(4.32)$$

\noindent a constant, independent of wealth and the ruin level.

As in Case A, we have five unknown constants with five equations that determine them.  The procedure for calculating these five constants and the resulting $h$ and $\pi^*$ on $\bf R$ is identical to the procedure described there, so we omit the details.

\medskip

\noindent{\bf C.} $\rho < r$.  In this case, $w^s = \bar c/(r - \rho)$; that is, the safe level is finite.  The general solution of (4.11) is given by (4.13), except that the positive root $\hat B_1$ is greater than 1 because $\rho < r$.  Also, $h(w) = 0$ for $w \ge w^s$.  The remainder of this case follows exactly as in Case A, except that $d = \hat B_1/(\hat B_1 - 1) > 1$ in (4.15)-(4.17) because $\hat B_1 > 1$.  It follows that the optimal investment in the risky asset decreases with respect to wealth on $(\kappa, w^s)$.

\medskip

We have the following proposition that reflects some of the work in this section.

\prop{4.1} {If the rate of consumption $c$ is given by $(4.1)$, then $(2.11)$ has a unique nonincreasing, convex solution $h \in C^1({\bf R}) \cap C^2({\bf R} - \{w^s\})$.}

\medskip

\noindent {\bf Remark:} Proposition 4.1 holds for more general consumption functions than the one given in (4.1).  Certainly, (2.11) has a solution for any piecewise linear rate of consumption by the same sort of argument used in this section.  For other types of consumption functions, the dual formulation will result in a nonlinear second-order differential equation with the nonlinearity in the first derivative.  Whether such a differential equation has a solution will depend upon the nature of the consumption function.  For example, if the nonlinear function satisfies a Lipschitz condition, then we can generally get a {\it local} solution, but we want, of course, a global solution; Tenenbaum and Pollard (1963, Theorem 62.22).

\medskip

In fact, we have shown more in this section than the existence of a solution to (2.11).

\th{4.2} {Suppose the rate of consumption $c$ is such that $(2.11)$ has a solution $h \in C^1({\bf R}) \cap C^2({\bf R} - \{w^s\})$.  If, in addition, $c(w) = \bar c + \rho w$ on $(\kappa, w^s),$ in which $\kappa$ is as specified in Theorem $3.1$ and $\rho > 0$, then the optimal investment strategy $\pi^*$ on $(\kappa, w^s)$ corresponding to $h$ is given by}

$$\pi^*(w) = \cases{{\mu - r \over \sigma^2}  {w + {\bar c \over \rho - r} \over 1 - d}, &if $\rho \ne r$, \cr
{\mu - r \over \sigma^2} {\bar c \over \delta + \lambda}, &if $\rho = r$,} \eqno(4.33)$$

\noindent {\it in which $d$ is given by $(4.16)$.}

\medskip

\pf Consider $\psi = \psi(w; \kappa)$, the minimum probability of ruin with ruin level $b = \kappa$.  The optimal investment strategy for $\psi$ is given by (4.33), as we showed in this section.  By Corollary 2.3, $\pi^*$ is (4.33) is also the optimal investment strategy for $h$ on $(\kappa, w^s)$.  $\square$

\subsect{4.2. Implied Utility Functions}

In this section, we show that the utility functions implied by equating $\pi$ and $c$ in the two optimization problems exhibit hyperbolic absolute risk aversion (HARA).  We calculate the implied utility functions by considering in detail the three cases: $\rho > r$, $\rho = r$, and $\rho < r$, as in Section 4.1.  Note that the domain $w \in (\kappa, w^s)$ corresponds to $c \in (\bar c + \rho \kappa, \bar c + \rho w^s)$ because $c(w) = \bar c + \kappa w$ on $(\kappa, w^s)$.

We first present and prove our main theorem of this paper, and then we remark on the results.  Recall that hyperbolic absolute risk aversion means that the absolute risk aversion is a hyperbolic function, in which the absolute risk aversion is given by $R_A(c) \equiv -u''(c)/u'(c)$ (Pratt, 1964).

\th{4.3}{Suppose that $(2.11)$ has a convex solution $h \in C^1({\bf R}) \cap C^2({\bf R} - \{w^s\}),$ with corresponding optimal investment strategy $\pi^*$ given in $(2.13)$.  Suppose that $(3.1)$ has a solution with corresponding optimal investment  and consumption strategies, $\pi^u$ and $c^u,$ respectively.  Let $\kappa,$ $w^s,$ and $w_\infty$ be as given in the hypothesis of Theorem 3.1.
\item{$(1)$} If $\pi^*(w) = \pi^u(w)$ and if $c(w) = c^u(w)$ on $(\kappa, w^s),$ then $c'(w) = \rho$ on $(\kappa, w^s)$ and $u$ exhibits hyperbolic absolute risk aversion on $(\bar c + \rho \kappa, \bar c + \rho w^s)$ for some constant $\bar c$.  Specifically, up to $(positive)$ linear transformation of $u,$ we have
\item\item{$(a)$} if $\rho > r,$ then
$$u(c) = {1 \over d} \left( c + {\bar c r \over \rho - r} \right)^d, \quad c \in (\bar c + \rho \kappa, \infty), \eqno(4.34)$$
\item\item{} \indent in which $d < 0$ is given by $(4.16);$
\item\item{$(b)$} if $\rho = r,$ then
$$u(c) = - {\bar c r \over \delta + \lambda} \exp \left( - {\delta + \lambda \over \bar c r} c \right), \quad c \in (\bar c + \rho \kappa, \infty); \eqno(4.35)$$
\item\item{$(c)$} if $\rho < r,$ then
$$u(c) = - {1 \over d} \left( {\bar c r \over r - \rho} - c \right)^d, \quad c \in (\bar c + \rho \kappa, r \bar c/(r - \rho)), \eqno(4.36)$$
\item\item{} \indent in which $d > 1$ is given by $(4.16)$.
\item{$(2)$} Conversely, if $u$ is as in items $(a), (b),$ or $(c)$ in $(1),$ and if $c(w) = c^u(w)$ on $(\kappa, w^s),$  then $\pi^*(w) = \pi^u(w)$ on $(\kappa, w^s)$.}

\medskip

\pf  This result follows easily from Theorems 3.1 and 4.2 by substituting $\pi^*$ as given in (4.33) into the expression for $u$ given by (3.10).  $\square$

\medskip

\hangindent 20 pt \noindent{\bf Remarks: (1)} If $\rho > r,$ then $u$ exhibits {\it decreasing} (hyberbolic) absolute risk aversion.  Indeed, the absolute risk aversion in this case equals

$$R_A(c) = {1 - d \over c + {\bar c r \over \rho - r}}, \eqno(4.37)$$

\item{} a decreasing function of consumption.  One implication of decreasing absolute risk aversion is that one is willing to spend less for insurance as wealth increases, an intuitively pleasing result (Pratt, 1964).  Also, the amount invested in the risky asset increases with wealth, as seen in (4.17).  Therefore, if we were to observe an individual who consumes at a linear rate that grows more quickly than the riskless rate, we could suppose that she is minimizing some function of her lifetime minimum wealth {\it or} maximizing her utility of lifetime consumption under a HARA utility function with decreasing $R_A$.

\item{} \indent If $\bar c = 0$, then $u$ in (4.34) is power utility, and $c(w) = \rho w$ on $(\kappa, \infty)$.  The relative risk aversion is defined by $R_R(w) = w R_A(w)$.  Indeed, the relative risk aversion in this case equals $1 - d$, which is greater than 1 because $d < 0$.

\medskip

\item{\bf (2)} If $\rho = r,$ then $u$ exhibits {\it constant} absolute risk aversion.  Also, note that $\bar c > 0$ in order for $u$ to be concave. One implication of constant absolute risk aversion is that the amount one is willing to spend for insurance is independent of wealth, a counterintuitive result (Pratt, 1964).  Also, the amount invested in the risky asset is constant, as seen in (4.32).  Such a utility function as given in (4.35) is not generally deemed acceptable in modeling preferences.  However, if we were to observe an individual who consumes at a linear rate that grows at the riskless rate, we could suppose that she is minimizing some function of her lifetime minimum wealth.

\medskip

\item{\bf (3)} If $\rho < r$, then $u$ exhibits {\it increasing} (hyberbolic) absolute risk aversion.  One implication of increasing absolute risk aversion is that one is willing to spend more for insurance as wealth increases, a counterintuitive result (Pratt, 1964).  Also, the amount invested in the risky asset decreases with wealth, as seen in (4.17).  For these reasons, such a utility function as given in (4.36) is not generally deemed acceptable in modeling preferences.  However, if we were to observe an individual who consumes at a linear rate that grows more slowly than the riskless rate, we could suppose that she is minimizing some function of her lifetime minimum wealth.  For such a problem with a finite safe level $w^s$, it is entirely reasonable that the individual would behave more conservatively as wealth increases towards the safe level.

\sect{5. Summary}

In this paper, we examined when the two problems of minimizing a nonnegative, nonincreasing function of lifetime minimum wealth and of maximizing utility of lifetime consumption result in the same optimal investment strategy.  We showed that the investment strategies are equal only when the consumption function is linear in wealth.  It followed that the corresponding investment strategy is also linear in wealth and the implied utility function exhibits hyperbolic absolute risk aversion (HARA).

Determining when the two investment strategies coincide is important because expected utility maximization is a well-established paradigm in the study of optimal investment and consumption strategies (especially under HARA utility).  Thus, obtaining the identical investment strategy under a different, potentially more intuitively appealing, criterion is helpful in validating the new criterion.  Also, as we observed in items (2) and (3) of the above Remark, if an individual invests as if her absolute risk aversion is nondecreasing, then one could interpret her behavior more naturally as if she were minimizing her lifetime probability of ruin (or more generally, a nonnegative, nonincreasing function of her lifetime minimum wealth) under an appropriate consumption function.

In the process of solving this problem, we proved a verification lemma for the minimum probability of ruin $\psi$ for a general consumption function, via an argument that avoids the usual verification argument for $L^\infty$ control problems that exhibit a Neumann condition.  We, then, proved that the optimal investment strategy is identical to the one when minimizing the expectation of a nonincreasing, nonnegative function $f$ of lifetime minimum wealth.  We used this correspondence to obtain a useful representation for $V^f$ in terms of $\psi$.  By using the Legendre transform, we explicitly calculated $\psi$ when the consumption is piecewise linear.

In summary, the two major contributions of this paper are (1) the optimal investment strategy corresponding to $V^f$ is independent of the function $f$, and (2) if the investors under the two problems considered in this paper behave similarly, then the utility function exhibits HARA.

\sect{Appendix:  Proof of Lemma 2.1}

For an arbitrary strategy $\pi \in {\cal A}$, let $W^\pi$ denote the wealth process when we use $\pi$ as the investment policy; similarly define $M^\pi$.  Define the hitting times $\tau_b \equiv \inf \{ t > 0: W^\pi_t = b \}$ and  $\tau_{w^s} \equiv \inf \{ t > 0: W_t^\pi = w^s \}$. (Technically, we should apply the superscript $\pi$ to these $\tau$'s, but we omit it because the notation is otherwise too cumbersome.)  If $W^\pi_t > b$ for all $t \ge 0$, then by convention $\tau_b = \infty$.

Because the time of death of the individual $\tau_d$ is independent of the Brownian motion, we can write $\psi$ as

$$\eqalign{\psi(w; b) &= \inf_{\pi \in {\cal A}} {\bf E}^w \int_0^\infty \lambda e^{-\lambda s} \, {\bf 1}_{\{M^\pi_s \le
b \}} \, ds \cr
& = \inf_{\pi \in {\cal A}} {\bf E}^w \int_{\tau_b}^\infty \lambda e^{-\lambda s} \, {\bf 1}_{\{M^\pi_s \le b\}} \, {\bf 1}_{\{ \tau_b < \tau_{w^s}\}}  \, ds \cr
&= \inf_{\pi \in {\cal A}} {\bf E}^w  \int_{\tau_b}^\infty \lambda e^{-\lambda s} \, {\bf 1}_{\{ \tau_b < \tau_{w^s}\}} \, ds = \inf_{\pi \in {\cal A}} {\bf E}^w \left[ e^{-\lambda \tau_b} {\bf 1}_{\{ \tau_b < \tau_{w^s}\}} \right].}  \eqno({\rm A}.1)$$

By using this formulation of the problem, the verification lemma follows from classical arguments, as we proceed to demonstrate.  First, define the stopping time $\tau_n \equiv \inf \{ t \ge 0: \int_0^t \pi^2_s \, ds \ge n \}$.  Then, define the stopping time $\tau^{(n)} = \tau_b \wedge \tau_{w^s} \wedge \tau_n$.

Assume that we have $v$ and $\beta$ as specified in the statement of Lemma 2.1.  Note that $v(w^s) = 0$ because if the individual reaches the safe level $w^s$, then her wealth cannot reach $b < w^s$ before she dies.  By applying It\^o's lemma to the function $f$ given by $f(w, t) = e^{-\lambda t} v(w)$, we have

$$\eqalign{& e^{-\lambda \left( t \wedge \tau^{(n)} \right)} v(W^\pi_{t \wedge \tau^{(n)}}) = v(w) - \lambda \int_0^{t \wedge \tau^{(n)}} e^{-\lambda s} \, v(W^\pi_s) ds + \int_0^{t \wedge \tau^{(n)}} e^{-\lambda s} \, v'(W^\pi_s)\sigma \pi_s dB_s \cr
&\qquad \qquad \qquad + \int_{0}^{t \wedge \tau^{(n)}} e^{-\lambda s} \left( (rW^\pi_s + (\mu-r) \pi_s - c(W^\pi_s))v'(W^\pi_s) + {1 \over 2} \sigma^2 \pi^2_s v''(W^\pi_s) \right) ds \cr
&\qquad = v(w) + \int_0^{t \wedge \tau^{(n)}} e^{-\lambda s} \, {\cal L}^{\pi_s}v(W^\pi_s) ds + \int_{0}^{t \wedge \tau^{(n)}} e^{-\lambda s} \, v'(W^\pi_s) \, \sigma \pi_s dB_s, \cr} \eqno({\rm A}.2)$$

\noindent in which the second equality follows from the definition of ${\cal L}^{\alpha}$ in (2.7).

If we take the expectation of both sides, the expectation of the last term in (A.2) is zero because

$${\bf E}^w \left[ \int_0^{t \wedge \tau^{(n)}} e^{-2 \lambda s} (v'(W^\pi_s))^2 \sigma^2 \pi_s^2 \, ds \right] \le (v'(b))^2 \sigma^2 \, {\bf E}^w \left[ \int_0^{t \wedge \tau^{(n)}} \pi_s^2 ds \right] < \infty, \eqno({\rm A}.3)$$

\noindent since $v'(w)$ is bounded by $| v'(b) |$ on $[b, w^s)$ because $v$ is decreasing and convex.  Thus, we have

$${\bf E}^w \left[e^{-\lambda (t \wedge \tau^{(n)})} v(W^{\pi}_{t \wedge \tau^{(n)}}) \right] = v(w)+ {\bf E}^w \left[\int_0^{t \wedge \tau^{(n)}} {\cal L}^{\pi_s} v(W^{\pi}_s) ds \right] \ge v(w), \eqno({\rm A}.4)$$

\noindent where the inequality follows from assumption (iii) of the lemma.

Because $v$ is bounded, $v(w^s) = 0$, and $v(b) = 1$, it follows from (A.4) and the Dominated Convergence Theorem (Royden, 1968) that

$$v(w) \le {\bf E}^w \left[e^{-\lambda \tau_b} v(W^{\pi}_{\tau_b}) {\bf 1}_{\{\tau_b < \tau_{w^s} \}} \right] = {\bf E}^w \left[e^{-\lambda \tau_b} {\bf 1}_{\{ \tau_b < \tau_{w^s} \}} \right], \eqno({\rm A}.5)$$

\noindent for any $\pi \in {\cal A}$.  Thus, it follows from (A.1) that $v \le \psi$.

Now, let $\beta$ be as specified in the statement of this lemma; that is, $\beta$ is the minimizer of ${\cal L}^\pi v$.   It follows from the above argument that we will have equality in (A.5), from which it follow that $v = \psi$.  Hence, we have demonstrated (2.8) and (2.9) on $[b, w^s)$.  Assumption (v) completes the proof.  $\square$


\bigskip
\centerline{\bf Acknowledgements} \medskip
This research of the first author is supported in part by the National Science Foundation under grant DMS-0604491.  We both thank Kristen Moore, David Promislow, and Thaleia Zariphopoulou for their insight.  We also thank the editor Martin Schweizer, an associate editor, and two anonymous referees for very helpful comments.

\vfill \eject

\sect{References}

\smallskip \noindent \hangindent 20 pt Al-Gwaiz, M. A.: Theory of distributions. Monographs and Textbooks in Pure and Applied Mathematics 159. New York: Marcel Dekker (1992)

\smallskip \noindent \hangindent 20 pt Barles, G., Daher, C., Romano, M.: Optimal control on the $L^\infty$ norm of a diffusion process. SIAM Journal on Control and Optimization {\bf 32}(3), 612-634 (1994)

\smallskip \noindent \hangindent 20 pt Bayraktar, E., Young, V. R.: Minimizing the lifetime shortfall or shortfall at death. Working Paper, Department of Mathematics, University of Michigan (2005)

\smallskip \noindent \hangindent 20 pt Bayraktar, E., Young, V. R.: Minimizing the probability of lifetime ruin under borrowing constraints. Insurance: Mathematics and Economics. To appear (2006)

\smallskip \noindent \hangindent 20 pt Browne, S.: Optimal investment policies for a firm with a random risk process: Exponential utility and minimizing the probability of ruin. Mathematics of Operations Research {\bf 20}(4), 937-958 (1995)

\smallskip \noindent \hangindent 20 pt Browne, S.: Survival and growth with a liability: Optimal portfolio strategies in continuous time. Mathematics of Operations Research {\bf 22}(2), 468-493 (1997)

\smallskip \noindent \hangindent 20 pt Browne, S.: Beating a moving target: Optimal portfolio strategies for outperforming a stochastic benchmark. Finance and Stochastics {\bf 3}, 275-294 (1999a)

\smallskip \noindent \hangindent 20 pt Browne, S.: The risk and rewards of minimizing shortfall probability, Journal of Portfolio Management {\bf 25}(4), 76-85 (1999b)

\smallskip \noindent \hangindent 20 pt Fleming, W. H., Soner, H. M.: Controlled Markov processes and viscosity solutions. New York: Spring-Verlag (1993)

\smallskip \noindent \hangindent 20 pt Fleming, W. H., Zariphopoulou, T.: An optimal investment/consumption model with borrowing. Mathematics of Operations Research {\bf 16}(4), 802-822 (1991)

\smallskip \noindent \hangindent 20 pt Heinricher, A. C., Stockbridge, R. H.: Optimal control of the running max. SIAM Journal on Control and Optimization {\bf 29}(4), 936-953 (1991)

\smallskip \noindent \hangindent 20 pt Hipp, C., Plum, M.: Optimal investment for insurers, Insurance: Mathematics and Economics {\bf 27}, 215-228 (2000)

\smallskip \noindent \hangindent 20 pt Hipp, C., Taksar, M.: Stochastic control for optimal new business. Insurance: Mathematics and Economics {\bf 26}, 185-192 (2000)


\smallskip \noindent \hangindent 20 pt Karatzas, I., Shreve, S.: Methods of mathematical finance. New York: Springer-Verlag (1998)

\smallskip \noindent \hangindent 20 pt Merton, R. C.: Continuous-time finance, revised edition. Cambridge, Massachusetts: Blackwell Publishers (1992)

\smallskip \noindent \hangindent 20 pt Milevsky, M. A., Ho, K., Robinson, C.: Asset allocation via the conditional first exit time or how to avoid outliving your money. Review of Quantitative Finance and Accounting {\bf 9}(1), 53-70 (1997)

\smallskip \noindent \hangindent 20 pt Milevsky, M. A., Moore, K. S., Young, V. R.: Asset allocation and annuity-purchase strategies to minimize the probability of financial ruin. Mathematical Finance {\bf 16}(4), 647-671 (2006)

\smallskip \noindent \hangindent 20 pt Milevsky, M. A., Robinson, C.: Self-annuitization and ruin in retirement, with discussion. North American Actuarial Journal {\bf 4}(4), 112-129 (2000)

\smallskip \noindent \hangindent 20 pt Pratt, J. W.: Risk aversion in the small and in the large. Econometrica {\bf 32}, 122-136 (1964)


\smallskip \noindent \hangindent 20 pt Royden, H. L.: Real analysis, second edition. New York: Macmillian (1968)

\smallskip \noindent \hangindent 20 pt Schmidli, H.: Optimal proportional reinsurance policies in a dynamic setting. Scandinavian Actuarial Journal {\bf 2001}(1), 55-68 (2001)

\smallskip \noindent \hangindent 20 pt Tenenbaum, M., Pollard, H.: Ordinary differential equations. New York: Dover Publications (1963)

\smallskip \noindent \hangindent 20 pt Young, V. R.: Optimal investment strategy to minimize the probability of lifetime ruin. North American Actuarial Journal {\bf 8}(4), 105-126 (2004)

\smallskip \noindent \hangindent 20 pt Zariphopoulou, T.: Consumption-investment models with constraints. SIAM Journal on Control and Optimization {\bf 32}(1), 59-85 (1994)

\bye